\documentclass{article}

\usepackage{arxiv}
\usepackage{amsmath}
\usepackage{algorithm}
\usepackage{algpseudocode}
\usepackage{adjustbox}
\usepackage{color}
\usepackage{graphicx,subfigure}
\usepackage{comment}
\usepackage{bm}

\usepackage{makecell}
\usepackage[utf8]{inputenc}
\usepackage[T1]{fontenc}   
\usepackage{hyperref}      
\usepackage{url}            
\usepackage{booktabs}       
\usepackage{amsfonts}       
\usepackage{nicefrac}      
\usepackage{microtype}      
\usepackage{lipsum}
\usepackage{graphicx}
\graphicspath{{media/}}   
\usepackage{geometry}
\title{Reduced Order Modeling of Partial Differential Equations on Parameter-Dependent Domains Using Deep Neural Networks}

\author{
  Martina Bukač  \\
Applied and Computational Mathematics and Statistics, University of Notre Dame, USA    \\
   \And
 Iva Mikuš \\
  University of Zagreb Faculty of Electrical Engineering and Computing, Croatia 
    \\
    \And
     Boris Muha \\
Department of Mathematics, Faculty of Science, University of Zagreb, Croatia    \\
    \And
    Domagoj Vlah \\
  University of Zagreb Faculty of Electrical Engineering and Computing, Croatia  
    \\
}

\begin{document}

\maketitle

\begin{abstract}
Partial differential equations (PDEs) are widely used for modeling various physical phenomena.  These equations often depend on certain parameters, necessitating either the identification of optimal parameters or the solution of the equations over multiple parameters. Performing an exhaustive search over the parameter space requires solving the PDE multiple times, which is generally impractical. To address this challenge, reduced order models (ROMs) are built using a set of precomputed solutions (snapshots) corresponding to different parameter values. These models allow for a fast approximation of the PDE solution when new parameters are introduced.

Recently, Deep Learning ROMs (DL-ROMs) have been introduced as a new method to obtain ROM, offering improved flexibility and performance. In many cases, the domain on which the PDE is defined also varies and may be described by parameters or measurements. Capturing this variation is important for building accurate ROMs but is often difficult, especially when the domain has a complex structure or changes topology.

In this paper, we propose a Deep-ROM framework that can automatically extract useful domain parametrization and incorporate it into the model. Unlike traditional domain parameterization methods, our approach does not require user-defined control points and can effectively handle domains with varying numbers of components.  It can also learn from domain data even when no mesh is available, which is often the case in biomedical applications. Using deep neural networks, our approach reduces the dimensionality of both the PDE solution and the domain representation, making it possible to approximate solutions efficiently across a wide range of domain shapes and parameter values.
We demonstrate that our approach produces parametrizations that yield solution accuracy comparable to models using exact parameters. Importantly, our model remains stable under moderate geometric variations in the domain, such as boundary deformations and noise—scenarios where traditional ROMs often require remeshing or manual adjustment.

\end{abstract}


\section{Introduction}
Applications such as design optimization and simulations in real-world settings often require solving partial differential equations (PDEs) for a wide range of parameter values \cite{legresley2000airfoil}, \cite{colciago2014comparisons}, \cite{raveh2001reduced}, \cite{aguado2015real}. Typically, given a parameter vector $\boldsymbol{\lambda}$, the goal is to find a solution $\boldsymbol{u}$ such that:
 \begin{align}
     \mathcal{L}(\boldsymbol{u};\boldsymbol{\lambda}) = f(\cdot; \boldsymbol{\lambda}) \text{ in $\Omega(\boldsymbol{\lambda})$}, \\
     \mathcal{B}(\boldsymbol{u};\boldsymbol{\lambda}) = g(\cdot; \boldsymbol{\lambda}) \text{ on $\partial \Omega(\boldsymbol{\lambda})$},
 \end{align}
 where $\mathcal{L}$ is a differential operator, $\mathcal{B}$ is a boundary operator, $f$ and $g$ are the source and boundary data, and $\Omega$ is a parameter-dependent domain. This task can be viewed as learning a solution map from the parameters $\boldsymbol{\lambda}$ to their corresponding solutions $\boldsymbol{u}$. A straightforward approach is to numerically solve the PDE for each parameter instance, but this becomes computationally prohibitive as the number of evaluations increases.

 To address this, reduced order models (ROMs) are commonly used to approximate the solution map with a much lower computational cost. These models are built on the assumption that solutions to the parametric PDE lie on a low-dimensional manifold \cite{quarteroni2015reduced}. The construction of such a ROM typically involves two phases. In the \textit{offline} phase, a dataset of high-fidelity parameter–solution pairs $(\boldsymbol{\lambda}_i, \boldsymbol{u}_i)$ is collected using methods such as the finite element method (FEM) \cite{quarteroni2008numerical} and ROM is obtained. In the \textit{online} phase, the ROM is evaluated on unseen parameters. Full-order solutions are then projected on a reduced basis, often obtained using techniques such as Proper Orthogonal Decomposition (POD) \cite{berkooz1993proper}. In the online phase, a reduced system is solved to obtain an approximate solution for a new parameter value, significantly reducing computation time.
 
 However, standard ROMs face limitations, especially when applied to PDEs defined on parameter-dependent domains. Constructing a projection-based ROM typically requires all solutions to be discretized over a common mesh or aligned computational domain, so that basis functions are defined consistently across parameter instances. This assumption breaks down when the domain itself varies with the parameters, which is common in many real-world scenarios.
 
 A typical workaround is to map each domain to a fixed reference configuration and transform the PDE accordingly \cite{benner2020model}. However, this mapping is often problem-specific and not easily generalized. Although this is feasible when domains vary smoothly and share the same topology, it becomes problematic when domains differ structurally, for instance, when holes appear or disappear, or when the number of connected components changes. In such cases, defining a single reference domain and constructing consistent mappings becomes difficult, or even impossible.
 
 To address these challenges, we propose a data-driven framework based on convolutional autoencoders to learn a compact domain representation directly from its characteristic function. We represent each domain as a binary image and use a convolutional autoencoder to extract a low-dimensional encoding. This encoding serves as a learned domain parametrization, which we then integrate into a neural reduced order model that predicts PDE solutions. This approach requires no meshing, no manual alignment, and no user-defined control points. It is naturally suited to problems where the domain is defined through measurements or where topological complexity varies across parameter instances. 
 Moreover, in many real-world scenarios—such as biomedical simulations, flow over irregular surfaces, or domains derived from imaging—constructing a computational mesh is often impractical or infeasible. This limits the applicability of traditional ROM techniques that rely on consistent meshing, graph connectivity, or control-point mappings to handle domain variability.

This challenge is particularly acute when domains are derived from imaging or experimental data, where mesh or graph structures may not be available at all. To overcome these limitations, we propose a fully data-driven, mesh-free approach that operates directly on binary images representing domain characteristic functions. Using convolutional autoencoders, our method learns compact geometric encodings without requiring explicit parametrization or meshing.

We demonstrate the effectiveness of this approach on synthetic two-dimensional PDE problems designed to reflect practical challenges, including geometric variability, boundary noise, and structural differences across domains.
 
\section{Related work}
  Recently, there has been an increasing interest in using Deep Learning for Reduced Order Modelling. Since neural networks are universal approximators \cite{cybenko1989approximation}, \cite{hornik1989multilayer}, \cite{daubechies2022nonlinear}, it follows that they can also approximate the solution map. In this section, various approaches that use deep learning for the construction of the ROM, as well as some classical approaches, are described. Throughout this section, it is assumed that the training set of parameter-solution pairs $(\boldsymbol{\lambda}_i, \boldsymbol{u}_i), i = 1, 2, \ldots N$ is obtained via measurements or a numerical method such as the Galerkin finite element method \cite{quarteroni2008numerical}. 
  \subsection{Reduced Order Model via Proper orthogonal decomposition }
  Proper orthogonal decomposition (POD) \cite{berkooz1993proper,chatterjee2000introduction} is a linear reduction technique. The basic assumption is that even though discretized solutions belong to a high-dimensional vector space of dimension $M$, solutions can be well approximated by a lower-dimensional vector space of dimension $k \ll M$. In this setting, all solutions can be vectorized and assembled into a matrix $\mathbf{U} \in \mathbf{\mathbb{R}^{M \times N}}$. To obtain the basis of such a lower-dimensional vector space, SVD decomposition  \cite{golub1971singular} is calculated, which means that orthogonal matrices $\mathbf{\hat{U}} \in \mathbf{\mathbb{R}^{M \times M}}, \mathbf{\hat{V}} \in \mathbf{\mathbb{R}^{N \times N}} $ and diagonal $\Sigma \in \mathbf{\mathbb{R}^{M \times N}}$ with non-negative diagonal elements $\sigma_1 \geq \sigma_2 \geq \ldots \geq \sigma_{min\{M,N\}} $, are found such that
  \begin{align*}
      \mathbf{U} = \mathbf{\hat{U}\Sigma \hat{V}}.
  \end{align*}
  The matrix of orthogonal basis vector for lower-dimensional vector space $V$ is then defined by extracting the first $k$ columns of matrix $\mathbf{\hat{U}}$. \\
  We can then assume that any solution of this parametric problem can be written as
  \begin{align*}
      \mathbf{u} = \mathbf{Vx},
  \end{align*}
where $\boldsymbol{x}$ is $k$-dimensional vector. Now, instead of finding $\mathbf{u}$, we have to find $\mathbf{x}$. Then, when  $\mathbf{Vx}$ is substituted instead of $\mathbf{u}$ in the Galerkin FEM formulation, a new system with fewer unknowns is obtained. In linear PDE, this results in decreasing the number of equations. In non-linear cases, this cannot be trivially done, so techniques such as the Discrete Empirical Interpolation Method (DEIM) \cite{chaturantabut2010nonlinear} are developed. Further difficulties in employing this ROM arise when the domain depends on parameters because all discretized solutions need to have the same dimension and each coordinates of all solutions have to correspond to the same points in the domain. To obtain ROM on parametric domains $\Omega(\boldsymbol{\lambda})$, equations are mapped to the reference domain $\hat{\Omega}$. Mapping $\mathcal{M}(\cdot; \boldsymbol{\lambda})$ maps the reference domain into $\Omega(\boldsymbol{\lambda})$, so the correspondence of points 
across different domains is established. To construct such a mapping, one can use free-form deformation (FFD) \cite{sederberg1986free}, \cite{salmoiraghi2016advances}, \cite{rozza2018advances}, radial basis function interpolation (RBF) \cite{manzoni2012model}, \cite{morris2008cfd}, inverse distance weighting interpolation (IDW) \cite{shepard1968two}, \cite{witteveen2009explicit}, \cite{forti2014efficient}, or
 surface registration using currents \cite{ye2024data}. Furthermore, domains should be parametrized as well. To get parametrization, a set of control points $\boldsymbol{x}_1, \boldsymbol{x}_2, \ldots, \boldsymbol{x}_{c}$ around the domain is selected, and displacements of these points from their corresponding positions on the reference domain are calculated. Control points on each domain have to match the same points in the reference domain through appropriate mappings to get meaningful parameters. Control points have to be chosen carefully so that they capture the shape characteristics well. This imposes difficulties when one wants to obtain a single ROM for situations where there are a varying number of components or holes.
 
POD is naturally combined with neural networks for parametric problems, where POD is used for data compression, and neural networks are used to predict POD coefficients from parameters \cite{HESTHAVEN201855}, \cite{WANG2019289}. In \cite{Selimefendigil01022014}, \cite{Selimefendigil15052015} POD is combined with Generalized Regression Neural Networks (GRNN) that predict POD coefficients from the parameters. Although successful, this approach suffers from the same limitations with respect to parametric domains. 
 
 \subsection{Physical Informed Neural Networks for Reduced Order Model}
In \cite{CHEN2021110666}, the authors use a different approach to construct the ROM. They construct the ROM via POD. The goal of the model is to predict coefficients in the reduced basis from parameter values using Physics-Informed Neural Networks (PINNs) \cite{RAISSI2019686}. The loss function is the sum of the mean-squared error between predicted and correct coefficients and the mean-squared error of the equation residual of the reduced system. If not enough high-fidelity solutions are available, the model can be trained to minimize the equation residual loss only.

\subsection{Deep ROM (DL-ROM)}
\begin{figure}[H]
     \centering
     \includegraphics[scale = 0.5]{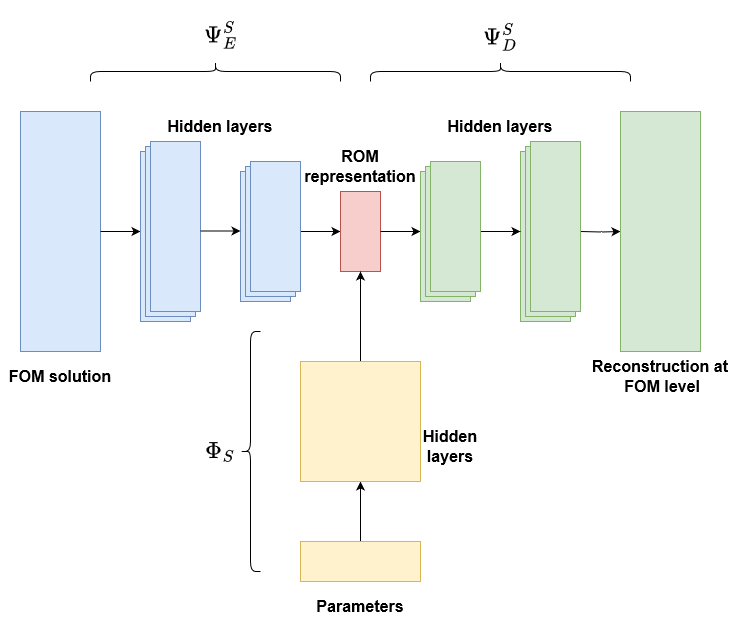}
     \caption{Basic DL-ROM architecture: encoder $\Psi_E^S$ and decoder $\Psi_D^S$ are trained to minimize reconstruction error between full order model solution and its reconstruction, After training, encodings of solutions are calculated with encoder $\Psi_E^S$ and Multilayer perceptron (MLP) $\Phi_S$ is trained to map parameters to encoding. To obtain a solution for a new parameter $\boldsymbol{\lambda}$, parameters are passed through  $\Phi_S$ and then through the decoder $\Psi_D^S$.}
     \label{basicROM}
 \end{figure}

 Franco et al. introduced Deep ROM (DL-ROM) in \cite{Franco_2022}. The idea is to obtain ROM by exploiting autoencoder architecture \cite{hinton2006reducing}. The model consists of three neural networks: an encoder $\Psi_E^S$,  a decoder  $\Psi_D^S$, and a multi-layer perceptron (MLP) \cite{haykin1998neural} $\Phi_S$. To obtain ROM, the encoder-decoder pair (i.e. autoencoder) is trained to minimize reconstruction error:
 \begin{align*}
     \frac{1}{N} \sum_{i=1}^N \| \boldsymbol{u}_i - \Psi_D^S(\Psi_E^S(\boldsymbol{u}_i))\|^2.
 \end{align*}
 The encoder for standard DL-ROM $\Psi_E^S(\boldsymbol{u}_i)$ maps the solution to a lower-dimensional space, called encoding, and the decoder lifts the encoding to the original solution size.  Finally, $\Phi_S$ is trained to minimize the mean square error between its output and encoder encoding:
 \begin{align*}
     \frac{1}{N} \sum_{i=1}^N\| \Phi_S( \boldsymbol{\lambda}_i) - \Psi_E^S(\boldsymbol{u}_i) \|^2. 
 \end{align*}
 The solution $\boldsymbol{u}$ for a new value of parameter $\boldsymbol{\lambda}$ is then:
 \begin{align}
     \boldsymbol{u} = \Psi_D^S( \Phi_S(\boldsymbol{\lambda})).
 \end{align}
 The entire architecture is shown on Figure \ref{basicROM}
 Notice that getting solutions from this kind of reduced-order model does not require solving any additional equations and can be efficiently done for multiple parameters at the same time, using GPU parallel computing. Authors used Convolutional Neural Networks (CNN) \cite{NIPS1989_53c3bce6} for the encoder and decoder and all domains involved were rectangular.
 \subsection{Deep Learning Proper Orthogonal Decomposition ROM}
Deep Learning Proper Orthogonal Decomposition ROM (POD-DL-ROM) was introduced in \cite{Fresca_2022}. The procedure goes as follows. Solutions $\boldsymbol{u}_i, i = 1, \ldots N$ are reshaped into vectors of length $N_h$ and assembled into a matrix $\mathbf{A}$ with $N_h$ rows and $N$ columns. Columns of matrix $\mathbf{A}$ are $\boldsymbol{a}_i, i = 1, 2, \ldots N$. The essence of the POD method is to get orthonormal vectors $\boldsymbol{\Phi_i}, i = 1, 2,  \ldots k$, $k \ll N$ that maximize:
 \begin{align*}
     \frac{1}{N} \sum_{i=1, \ldots N,j=1, \ldots k} | \langle \boldsymbol{a}_i, \boldsymbol{\Phi}_j\ \rangle |^2.
 \end{align*}
 Then, the following approximation hold:
 \begin{align*}
     \boldsymbol{u} \approx \mathbf{V}_k\boldsymbol{u}_k
 \end{align*}
Matrix $\mathbf{V}_k$ is usually found using SVD decomposition or randomized SVD (rSVD) \cite{rSVD}.  
The solution $u$ can be approximated with POD basis and represented as $k$-dimensional vector $[ \langle \boldsymbol{u}, \boldsymbol{\Phi}_1 \rangle , \ldots \langle \boldsymbol{u},  \boldsymbol{\Phi}_k \rangle]$. Then, similar to the black-box DL-ROM, the model consists of three networks: encoder $\Psi_E$, decoder $\Psi_D$, and parameter-encoding network $\Phi$. The difference is that the encoder and decoder are trained on POD coefficients, thus reducing the dimension of the problem. For this model to work on multiple domains, POD coefficients at the reference domain and mapping between domain and reference domain are needed.

\subsection{Graph convolutional autoencoder approach for DL-ROMs}
Following DL-ROM framework, in \cite{pichi2023graph} Graph convolutional autoencoders (GNNs) \cite{zhou2020graph}, \cite{monti2017geometric}, \cite{park2019symmetric}, \cite{zhou2020fully} are used for dimension reduction. Using the proposed autoencoder technique, one can obtain ROM for non-rectangle domains using the GNN autoencoder. However, graph pooling and a fully-connected layer work only when all domains have the same number of nodes and domain parametrization is known as well. In applications, especially in biomedicine, both constraints may not hold. 

In \cite{zhao2025diffeomorphism}, a harmonic mapping is used to map the domain and solution to a reference configuration. The neural operator is learned on a reference configuration and mapped to the original domain via inverse of the harmonic mapping. This approach may not be feasible when such harmonic mapping is difficult or impossible to construct due to a wide geometric or topological variability. \\
Similarly as in our work, in \cite{eichinger2020surrogate} a characteristic function of the domain is used, and the solutions are interpolated onto a rectangular grid using CNN. The input is a domain characteristic function and the output is the solution on that domain. Contrary to our work, this model is not capable of handling non-geometric parameters simultaneously.

\section{Our method}
In this work, we develop an efficient Reduced Order Model (ROM) that operates effectively across different domains, based on a Deep Learning Reduced Order Model (DL-ROM) approach. We refer to our method as Deep Learning Domain Aware ROM (DL-DA-ROM). A key part of this approach is learning a meaningful parametrization of the domain. In many problems, the domain may be defined through parameters, evolve with the system (e.g., moving domains), or come from measurements.

\begin{figure}[H]
     \centering

     \includegraphics[scale = 0.5]{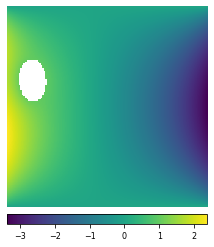}
     \includegraphics[scale = 0.5]{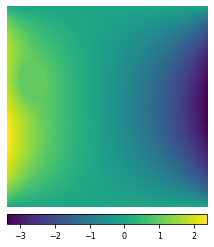}
     \caption{Left panel shows a sample function, and the right panel shows its interpolation into an equidistant rectangular grid.}
     \label{trainingSample}
 \end{figure}

Since the ROM needs to predict both the domain and the solution on that domain, we separate those two tasks.
We assume that the solution is at least continuous and can be interpolated onto rectangular grid without significant loss of accuracy, see Figure \ref{trainingSample}. This interpolation allows us to use CNNs to build an autoencoder for the solutions. This setup is especially useful when domains come from measurement data and creating a mesh is either difficult or too costly.
While it is possible to interpolate the solution and then mask out the region outside the domain, we found that using the full interpolated solution yields better performance. Interpolation naturally smooths the data, reducing sharp discontinuities introduced by irregular meshes or complex geometries. This smoothing effect makes the solution easier to approximate with CNNs, which are particularly well-suited for capturing smooth spatial patterns. As a result, CNN-based autoencoders can more efficiently learn a compact representation of the solution field without being disrupted by artificial discontinuities at domain boundaries.

\begin{figure}[H]
    \centering
    \includegraphics[scale = 0.2]{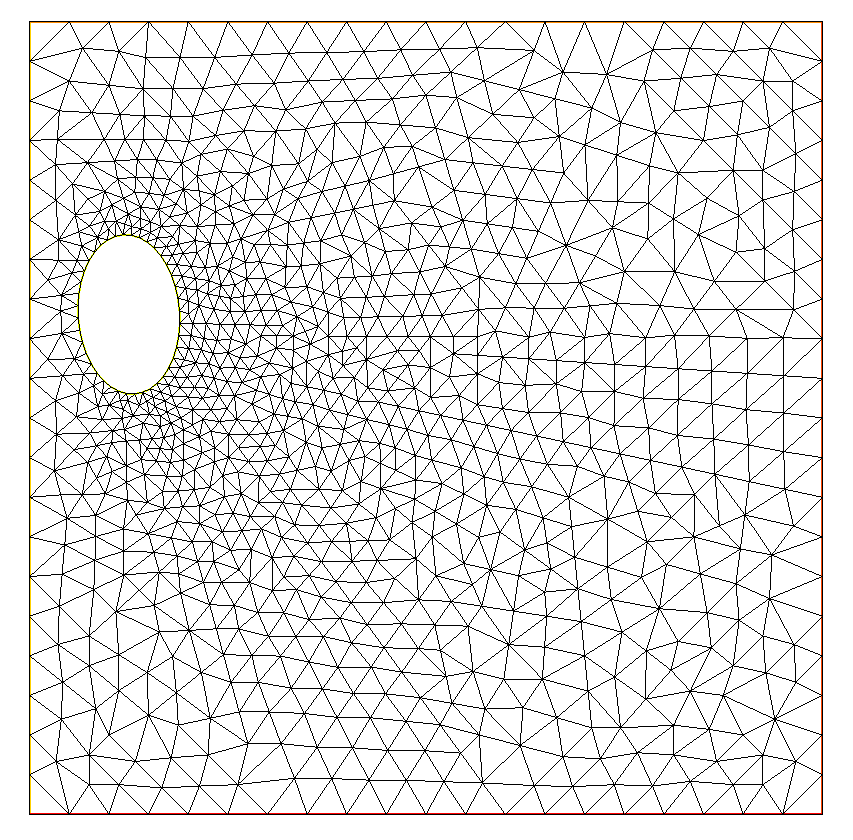}
    \includegraphics[scale = 0.55]{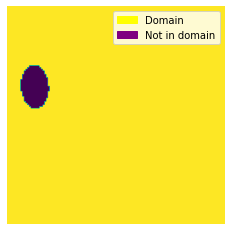}

    \caption{Comparison between a mesh-based domain representation (left) and the corresponding interpolated characteristic function represented as an image used by our method (right). The binary image representation allows us to learn a low-dimensional domain encoding without requiring a mesh.}
    \label{meshes}
\end{figure}

When domains are parametrically described, it is not necessary to assume that each domain has the same number of parameters; they may have varying numbers of subdomains or holes. Our approach assumes that each domain characteristic function can be interpolated onto a fixed rectangular grid, see Figure \ref{meshes}.  \\
To represent the domain, we use its characteristic function interpolated onto a rectangular grid, resulting in a binary image. This image is stored as a bitmap, where each pixel is assigned a value of 1 if the corresponding point lies inside the domain and 0 otherwise. A convolutional autoencoder is then trained on these bitmap images to extract a compact latent representation, which serves as a learned, low-dimensional domain parameterization. This mesh-free approach accommodates complex geometries and topological variations without requiring control points, meshing, or predefined mappings. To extract the domain parametrization, we train a convolutional autoencoder on these binary images, using the binary cross-entropy loss\cite{zhang2018generalized} to classify each pixel as part of the domain or not which is given by 
\begin{align*}
    l(x,y) = \frac{1}{nm} \sum_{i=1}^n 
    \sum_{j=1}^m -y_{ij}\log(x_{ij})-(1-y_{ij})\log(1-x_{ij}),
\end{align*}
where $\boldsymbol{y}$ is target bitmap image with $n$ rows and $m$ columns and $\boldsymbol{x}$ is network output image, whose pixels values $x_{ij}$ denote probability that pixel $x_{ij}$ is equal to one. The encoded domains are then used as the parametrization of the domain.
The low-dimensional latent code produced by the encoder is used as domain parametrization in our model. The autoencoder architecture is inspired by the encoder-decoder networks commonly used in semantic segmentation \cite{badrinarayanan2017segnet, noh2015learning}, with two classes that distinguish whether each pixel lies within or outside the domain. This approach can be extended to problems with multiple domains, such as fluid-structure interaction problems, providing simultaneous geometric parametrization for the entire system.

If domains are synthetically generated, it is also possible to train the domain autoencoder with a larger dataset than the solution autoencoder, which can help improve the quality of the learned parametrization and its generalization to unseen geometries.

\begin{figure}[H]
     \centering
     \includegraphics[scale=0.3]{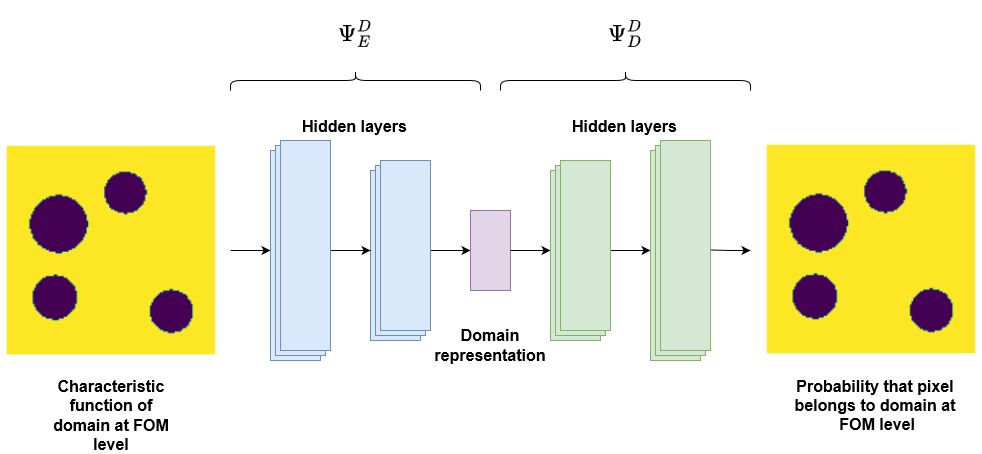}
     \caption{Domain Autoencoder used for domain parametrization. The input to the convolutional neural network is a characteristic function represented as a bitmap image. The output of the autoencoder is a rectangular grid where each pixel value represents the predicted probability that the pixel belongs to the domain. The low-dimensional encoder output (latent representation) is used as the domain parametrization.}
     \label{domainAE}
 \end{figure}

 The domain autoencoder is shown in Figure \ref{domainAE}.

\begin{figure}[H]
     \centering
     \includegraphics[scale=0.45]{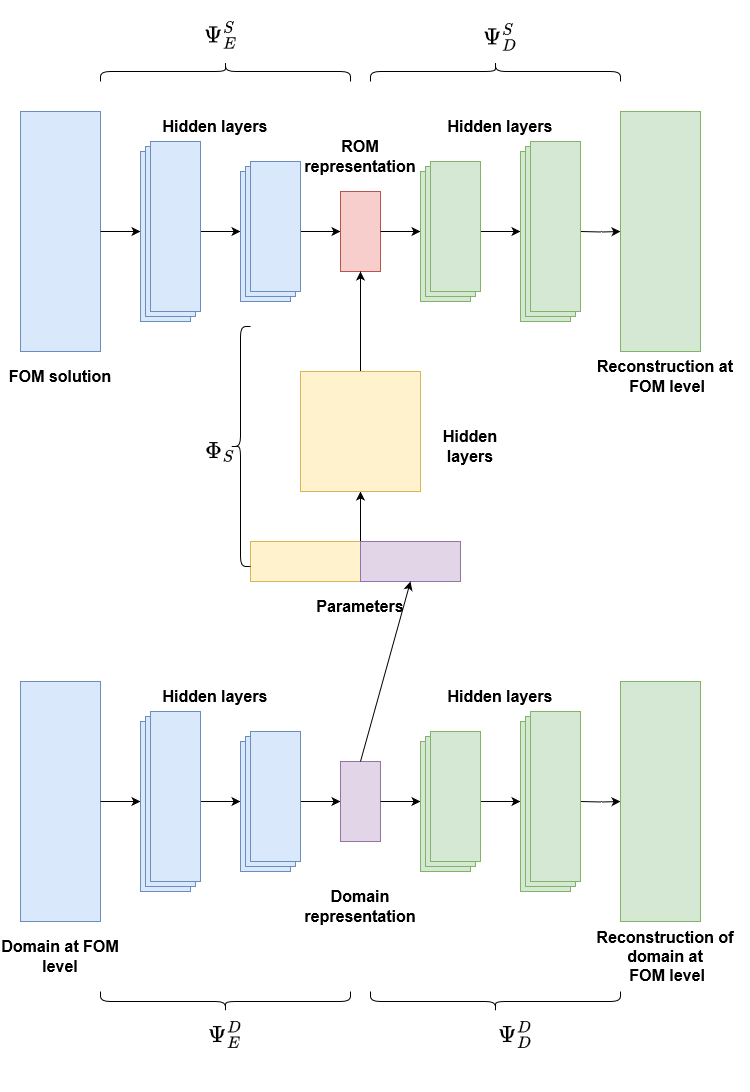}
     \caption{Deep ROM architecture for PDEs with parameter-dependent domains. In addition to the standard Deep ROM structure, our method incorporates a domain autoencoder (bottom) that encodes the characteristic function of the domain into a low-dimensional latent representation. This representation is combined with other problem parameters and passed to a fully connected network $\Phi_S$, which predicts the solution encoding. The predicted latent solution is then decoded using $\Psi^S_D$ to reconstruct the solution on the full-order model (FOM) grid. Both the solution and domain autoencoders are trained using CNNs. The domain representation serves as a learned geometric parametrization that allows the ROM to generalize across different domain shapes and configurations.}
     \label{domainROM}
 \end{figure}

 Fully-connected neural network $\Phi_S$ is trained to map input parameters (e.g., domain encoding and physical parameters) to the solution representation.  The predicted latent solution is then passed through the decoder $\Psi^S_D$ to reconstruct the full solution on the interpolated grid. The model is shown in Figure \ref{domainROM}  and all important terms are given in Table \ref{tableLabels}.\\

  \begin{center} 
\begin{table}[ht]
 \begin{adjustbox}{ center}

\begin{tabular}{ |c |c |}
 \hline
 \textbf{label} & \textbf{meaning} \\
 \hline 
   $\Omega(\lambda)$ & parameter-dependent domain \\ 
   \hline
   $\Psi_E^S$ & solution encoder \\
   \hline
   $\Psi_D^S$ & solution decoder \\
   \hline
  
   $\Phi_S$ & neural network that maps parameters into solution encoding \\
   \hline
    
   $\Psi_E^D$ & domain encoder \\
   \hline
    
   $\Psi_D^D$ & domain decoder \\   
   \hline

\end{tabular}
\end{adjustbox}

\caption{Notation used in this work. }
    \label{tableLabels}

\end{table}
\end{center}
 
Constructing the ROM-also known as the \textit{offline} phase, involves gathering the data, either using numerical simulation, measurements, or possibly a combination of both. During this phase, both the domain and solution autoencoders are trained to obtain meaningful low-dimensional representations. 
Before training the solution autoencoder training solution data is standardized by applying the following transformation:
\begin{align*}
    \boldsymbol{u}_i = \frac{\boldsymbol{u}_i- \mu_{s}}{\sigma_{s}},1, \ldots N,
\end{align*}
where $N$ is the number of samples in the training dataset and $\mu_{s}$ and $\sigma_{s}$ are the mean and standard deviation of the training dataset.
After the solution and the domain autoencoders are trained, we evaluate the domain encoder to extract the low-dimensional domain parameters. These are then concatenated with the PDE parameters. 
Prior to training the solution prediction network $\Phi_S$, the concatenated input equation and domain parameters are standardized:
\begin{align*}
    f_{ij} = \frac{f_{ij}- \mu_{j}^p}{\sigma_{j}^p},i=1,\ldots N,j = 1 \ldots k,
\end{align*}
where $f_{i,j}$ the $j$-th component of the concatenated parameter vector for the 
$i$-th training sample, $k$ is dimension of a parameter vector, $N$ number of samples and $\mu_{j}^p, \sigma_{j}^p$ are mean and standard deviation of the $j$-th  component of parameter vectors in the training set.  
The network $\Phi_S$ is then trained to map from this combined, standardized parameter vector to the corresponding solution encoding. The full procedure is summarized in Algorithm \ref{offline}. 
Once the training is complete, solutions for unseen parameters are calculated in the \textit{online} phase, as detailed in Algorithm \ref{online}.  Given a new domain and PDE parameters, the domain is first passed through the trained encoder $\Psi_E^D$, the standardized result is concatenated with the standardized PDE parameters, and the solution is reconstructed using the trained networks. The predicted ROM solution $\hat{u}$ is then transformed back to the original scale using:
\begin{align*}
    \boldsymbol{\hat{u}} = \sigma_s  \boldsymbol{\hat{u}} + \mu_s.
\end{align*}

\begin{algorithm}[H]
    \caption{Offline phase}\label{offline}
        \begin{algorithmic}
       \Require Parameter-solution-characteristic function triplets $(\boldsymbol{\lambda}_i, \boldsymbol{u}_i, \boldsymbol{c}_i),i = 1, \ldots N$, hyperparameters for training of networks $\Psi_E^S, \Psi_D^S$, domain autoencoder $\Psi_E^D, \Psi_D^D$, and parameter-encoding network $\Phi_S$
        \Ensure Trained networks
        \State Standardize solution data
        \State Train solution encoder $\Psi_E^S$ and decoder $\Psi_D^S$ to minimize reconstruction loss $\frac{1}{N} \sum_{i=1}^N\| \boldsymbol{u}_i - \Psi_D^S(\Psi_E^S(\boldsymbol{u}_i))\|^2$
        \State Evaluate trained $\Psi_E^S$ on $\boldsymbol{u}_i, i=1, \ldots N$ and get encodings $\boldsymbol{e}_i=\Psi_E^S(\boldsymbol{u}_i), i=1,\ldots N$ 
        \State Train domain encoder $\Psi_E^D$ and decoder $\Psi_D^D$ to minimize binary cross-entropy loss  
        \State Evaluate trained $\Psi_E^D$ and get domain parameters $\boldsymbol{d}_i = \Psi_E^D(\boldsymbol{c}_i), i=1,\ldots N$  
        \State Concatenate domain and equation parameters $\boldsymbol{d}_i$ and $\boldsymbol{\lambda}_i$ into new parameters $\boldsymbol{f}_i, i= 1, \ldots N$
        \State Standardize new parameters $\boldsymbol{f}_i, i = 1, \ldots N$
        \State Train $\Phi_S$ to minimize reconstruction loss $\frac{1}{N}\sum_{i=1}^N\|\Phi_S(\boldsymbol{f}_i) - \boldsymbol{e}_i\|^2$
        \State Save all trained networks
        \State Save mean $\boldsymbol{\mu}_s$ and standard deviation $\boldsymbol{\sigma_s}$ of training dataset and mean $\boldsymbol{\mu}_p$ and standard deviation $\boldsymbol{\sigma}_p$ vectors of training feature dataset
    \end{algorithmic} 
\end{algorithm}

\begin{algorithm}[H]
    \caption{Online phase}\label{online}
        \begin{algorithmic}
       \Require: Neural networks $\Psi_D^S, \Psi_E^D$, and $\Phi_S$, values of parameters and/or domains $\boldsymbol{d}$ for which we want to evaluate solution $(\boldsymbol{\lambda}_i, \boldsymbol{c}_i), i = 1, \ldots M $, mean and standard deviation of training parameters $\boldsymbol{f}_i$ $\boldsymbol{\mu}_p, \boldsymbol{\sigma}_p$, mean and standard deviation of solutions $\mu_s, \sigma_s$
        \Ensure Solutions for given parameters and domains $\boldsymbol{u}_i, i = 1,\ldots M$
        \State Compute domain parameters $\boldsymbol{d}_i = \Psi_E^D(\boldsymbol{c}_i)$
        \State Concatenate $\boldsymbol{f}_i = (\boldsymbol{\lambda}_i, \boldsymbol{d}_i), i= 1, \ldots M$
        \State Standardize $\boldsymbol{f}_i$ based on training mean and standard deviation
        \State Calculate ROM prediction $ \boldsymbol{\hat{u}} = \Psi_D^S(\Phi_S(\boldsymbol{f}_i))$
        \State Return solution $(\sigma_s \boldsymbol{\hat{u}} + \mu_s) \odot d$
     \end{algorithmic} 
\end{algorithm}

\section{Experiments and Results}
In this section we evaluate DL-DA-ROM performance on multiple problems, varying in domain parametrization complexity, domain topology, and PDE type. We compare performance against models using exact geometric parameters when such parametrizations are available and assess robustness to domain variability. We include a setting with variable domain topology, where classical fixed-dimensional parametrizations are not directly applicable, as well as cases where constructing mappings to a reference domain may be challenging in practice. We compare our approach with POD-DL-ROM. In all mentioned cases, constructing a classical ROM with diffeomorphic mapping to a reference domain is either infeasible due to variable domain topology, or impractical due to the complexity of the geometry deformations involved, particularly due to the requirement that the mapping maintains a positive Jacobian. Morover, if valid mapping to a reference domain is constructed, change of variables would result in the system of equations depending nolinearly on parameters, thus necessitating further approximations such as DEIM. 

Throughout this section, we will use relative error to evaluate our approach. The relative error $\epsilon$ for a sample $(\boldsymbol{\lambda}, \boldsymbol{u}_h, \boldsymbol{d}_h)$ is:

\begin{align*}
    \epsilon(\boldsymbol{u}_h, \hat{\boldsymbol{u}_h}) = \frac{\sqrt{\sum_{i,j}(u^{i,j}_h - \hat{u}^{i,j}_h)^2d^{i,j}_h}}{ \sqrt{\sum_{i,j} (u^{i,j}_h)^2d^{i,j}_h}},
\end{align*}
where $\boldsymbol{\lambda}$ is the parameter vector, $\boldsymbol{u}_h$ solution discretized on a mesh with $n \times m$ points, $\boldsymbol{\hat{u_h}}$ is DL-ROM prediction and $\boldsymbol{d}_h$ is an approximation of characteristic functions of domains on that same mesh. \\ 
For visualization purposes, we also define a pointwise relative error at position $(i,j)$:
\begin{align*}
    \epsilon_{i,j}(u_k, \hat{u}_h) = \frac{(u_h^{i,j} - \hat{u}^{i,j}_h)d^{i,j}_h}{\sqrt{\sum_{k,l} (u^{k,l}_h)^2d^{k,l}_h}},
\end{align*}
where $u^{i,j}_h, \hat{u}^{i,j}_h, d^{i,j}_h$ are discretized solutions, model approximation of the solution, and characteristic function value in \textit{i}-th row and \textit{j}-th column.

We use a FreeFEM++ software \cite{MR3043640} to generate finite element solutions and PyTorch \cite{NEURIPS2019_9015} to implement and train all neural networks. All solutions are projected onto a uniform rectangular grid to enable the use of CNNs. For the advection–diffusion examples, the elliptical-hole and circular-hole cases are discretized on a $128 \times 128$ grid, while the Navier–Stokes case uses a $256 \times 512$ grid. Interpolation relative error is calculated as L2 relative error of the data on interpolated mesh projected to original mesh and original data and is approximately $0.06\%$ for every example.
In all experiments, we use Leaky ReLU activations in convolutional layers and apply batch normalization \cite{ioffe2015batch} after each layer, except for the final layers of the encoder and decoder. The decoder mirrors the encoder in architecture and number of channels. Upsampling layers are used in selected layers to ensure output resolutions match those of the corresponding encoder layers. All networks are trained using Adam optimizer \cite{kingma2014adam} with the One Cycle Policy learning rate scheduler \cite{smith2019super}, with a maximum learning rate of $1 \times 10^{-3}$ for autoencoder networks.

All training, test, and validation sets used have $1000$ samples. 
The batch size for the autoencoder is $50$.
Autoencoder network hyperparameters are chosen by trial and error with respect to autoencoder reconstruction, where we searched on latent dimensions $15$, and $30$ and the number of convolutional channels for each encoder layer (the decoder mirrors the encoder). The hyperparameters tested are shown in Table \ref{AEgridSearchParams}, are identical across all examples. Autoencoders with more channels per layer have better training and test accuracy, so we picked the best two for further analysis.
\begin{center} 
\begin{table}[H] \centering
\begin{tabular}{ |c |c |}
 \hline
  
   Latent dimension & $15, 30$ \\ \hline
    \thead{Nbr. of channels in each  \\ encoder layer} & \thead{ $[4,4,8,8,16,16,16]$, \\   $[4 ,8,  16, 32, 32, 64, 64]$ , \\   $[8,16,16,32,64,64,64]$} \\  \hline 
    \thead{Strides} & \thead{$[1, 1,  2, 2, 2, 1, 1]$}  \\ \hline
    Maximum learning rate & 1e-3 \\ \hline
    Batch size & 50 \\  \hline
    Kernel size & 5 \\
 \hline 
\end{tabular}
\caption{ Autoencoder grid search hyperparameters. We trained 6 models with varying channel sizes and latent dimension sizes. The number of channels in decoder layers mirrors the encoder number of channels per layer. We have chosen two autoencoders with the smallest mean relative reconstruction errors on the validation set, both with $8,16,16,32,64,64,64$ channels in encoders and dimensions $15$ and $30$ for further analysis.}
    \label{AEgridSearchParams}
\end{table}
\end{center}
 Since $\Phi_S$ is usually small compared to autoencoder and fast to train, we employ grid search to choose the best hyperparameters for $\Phi_S$. Hyperparameters are shown in Table \ref{gridSearchParams}. The grid search was done for two ROM representations, obtained from solution autoencoders with $8,16,16,32,64,64,64$ channels per layer and latent dimensions equal to $15$ and $30$, respectively. Domain features are obtained from autoencoders trained on $9000$ domains with $ 8, 16, 16, 32, 64, 64, 64$ channels per layer in each layer and strides $1, 1, 2, 2, 2, 1, 1$ per layer, respectively. For each model, the One Cycle Policy scheduler with Adam optimizer was used.  We trained two different domain autoencoders, one with latent dimension $20$ and one with latent dimension $30$.  We choose the best hyperparameters with respect to the mean-squared error between the prediction of $\Phi_S$ and ROM representation on validation set. 
 \begin{center} 
\begin{table}[H] \centering
\begin{tabular}{ |c |c |}
 \hline
  
    Hyperparameter  & Values \\ \hline
  Training set size & 1000 \\ \hline 
  Number of epochs & 1500 \\ \hline 
  Number of layers & 1,2,3,4 \\ \hline 
  Neurons per layer & 64, 128, 256, 512, 1024, 2048 \\ \hline
  Dropout rate & 0, 0.025, 0.05, 0.1, 0.2 \\ \hline
  Maximum learning rate & 1e-5, 3.16e-5, 1e-4, 3.16e-4, 1e-3, 3.16e-3, 1e-2 \\ \hline
  Batch size & $8, 16, 32, 64, 128$ \\ 
 \hline 
\end{tabular}
\caption{Hyperparameters  grid search for $\Phi$ is performed, which results in $4200$ models in total.   }
    \label{gridSearchParams}
\end{table}
\end{center}

\subsection{Advection-diffusion on different domains}
In this section, we evaluate the accuracy and generalization capabilities of our proposed DL-DA-ROM method on a parameterized advection-diffusion problem. The goal is to compare the performance of models using exact domain parameters (e.g., ellipse descriptors) versus those using domain encodings learned by the autoencoder. This setting allows for a controlled comparison, as both the physical and geometric parameters are known. By analyzing prediction accuracy across different latent dimensions, we assess the quality of the learned domain representations and their ability to replace traditional parameterizations. Precisely, we consider the following advection-diffusion equation on parametric domains:
\begin{equation}
    \begin{cases}
      -\Delta u + 0.1 \mathbf{c} \cdot \nabla u = 1 & \text{on $\langle0, 1\rangle^2$ without ellipse}\\
      u = 1 & \text{ on edge of ellipse} \\
      u = 0 & \text {for $y=0$ and $y=1$} \\
      \frac{\partial u}{\partial \nu} = \beta & \text{ for $x=0$} \\
      \frac{\partial u}{\partial \nu} = -\beta & \text{ for $x=1$}      
    \end{cases}       
\end{equation}
Parameters are:

\begin{itemize}
    \item  transport angle $\phi \in   [0, 2 \pi]$ which determines vector $\mathbf{c}=(cos(\phi), sin(\phi))$ 
    \item  center of the ellipse $(x_0, y_0) \in [0.15, 0.85]^2$
    \item angle between major axis of the ellipse and the x-axis $\alpha \in [0, \pi]$ 
    \item boundary flux $\beta \in [1, 10]$
    \item half-axes $a, b \in [0.05, 0.25]$
\end{itemize}
In total, there are $7$ parameters: two \textit{PDE parameters}—transport angle $\phi$ and boundary flux $\beta$—and five geometric parameters defining the ellipse: center coordinates $(x_0, y_0)$, half-axes $a$ and $b$, and orientation angle $\alpha$. We assume that the ellipse has to be at least 0.01 units away from each of the boundaries to avoid geometric artifacts.

In this case, it is not feasible to construct a diffeomorphism from the physical domain to a fixed reference domain, for example to the square with a circular hole centered at $(0.5,0.5)$ and radius $0.15$, where the center is placed at the midpoint of the square and the radius is chosen as an average radius in the parameter distribution. More precisely, determinant of the Jacobian changes sign in the cases where ellipse approaches the square boundaries. While a diffeomorphism theoretically exists between any two such domains sharing the same topology, constructing it numerically via standard mesh morphing techniques (Laplacian, biharmonic, or elasticity-based smoothers) produces degenerate Jacobians for these parameter combinations, making the approach impractical for the full parameter range.

Since the domain is generated from known ellipse parameters, we can directly compare the performance of our model using either exact geometric parameters or those inferred by the domain autoencoder. For each type of domain parametrization, we train a fully connected neural network $\Phi_S$ with hyperparameters selected via grid search.

\begin{figure}[H]
     \centering
     \includegraphics[scale = 0.2 ]{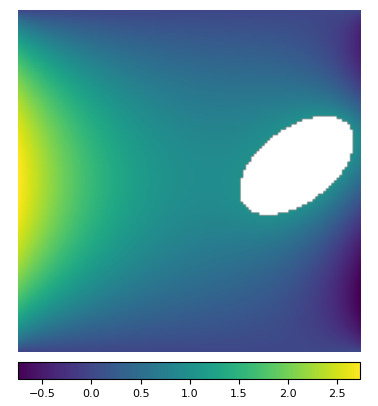} 
     \includegraphics[scale = 0.2]{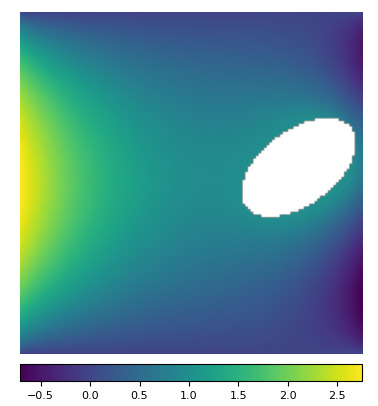} 
     \includegraphics[scale=0.2]{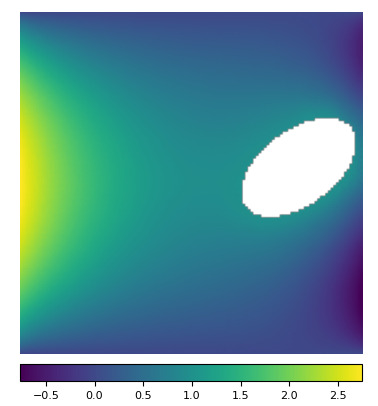}  
    \includegraphics[scale = 0.2]{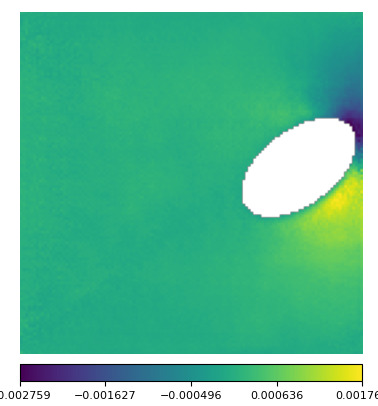}
    \includegraphics[scale = 0.2]{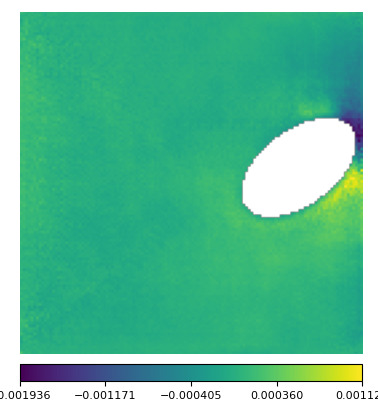}
 
     \caption{Comparison of predictions and relative errors for the DL-ROM with ROM encoding dimension $30$ and domain encoding dimension $20$. From left to right: Exact solution, network prediction for a model trained on exact domain parameters, network prediction for a model trained on approximated domain parameters, relative error for the model trained on exact domain parameters, and relative error for the model trained on approximated domain parameters. Results correspond to the ellipse with center $(0.815,\, 0.456)$, half-axes $a = 0.194$, $b = 0.110$, transport angle $\phi = 2.497$, and boundary flux $\beta = 6.884$. All models were trained on $1000$ samples.
 }
     \label{ellipse}
 \end{figure}

A comparison of the vanilla DL-ROM and our approach is shown in Table \ref{table1}. We compare DL-ROM performance when exact ellipse parameters are used with domain parametrization calculated with the autoencoder. 
Performance of the best models is shown in Figure \ref{ellipse}, evaluated on one of the examples from the test set.
 
\begin{center}
\begin{table}[H] 
\begin{adjustbox}{   center}
\begin{tabular}{ |c |c |c |c|c| }
 \hline
  
 & \thead{test error, \\ ROM encoding \\dimension 15, \\ domain encoding \\ dimension 20} 
 & \thead{test error, \\ ROM encoding \\ dimension 15, \\ domain  encoding \\ dimension 30} 
  &  \thead{test error, \\ ROM encoding \\ dimension 30, \\ domain encoding \\ dimension 20}
& \thead{test error, \\ ROM encoding \\ dimension 30, \\  domain encoding \\ dimension 30}   
\\
\hline
  \thead{exact ellipse \\ parameters only} & 0.0348
            &  0.0348   &  0.0383 &   0.0383                 
                             \\
   \hline
    
    \thead{calculated \\
    ellipse parameters} &  0.0339               
                      & 0.0361       &    0.0310    & 0.03975  \\
    \hline
\end{tabular}
\end{adjustbox}
    \caption{Comparison of DL-ROM and DL-DA ROM performance: Difference of the performance on relative test errors depending on the ROM encoding dimension and domain encoding dimension. We differentiate between cases where geometric parameters are exactly known (1st row) and when we use only geometric parameters obtained by autoencoder (2nd row). We can see that using the autoencoder domain parametrization does not degrade model performance. We report on the best model for each ROM encoding dimension and domain encoding dimension. Models trained on calculated domain parameters perform similar as models trained using only exact domain parameters. Even though results are similar, larger fully-connected neural network is needed in case of calculated domain parameters (4 layers with 2048 neurons vs 4 layers with 512 neurons).
  }
    \label{table1}
\end{table}
\end{center}

\begin{figure}[!h]
     \centering
     \includegraphics[width=0.4\linewidth]{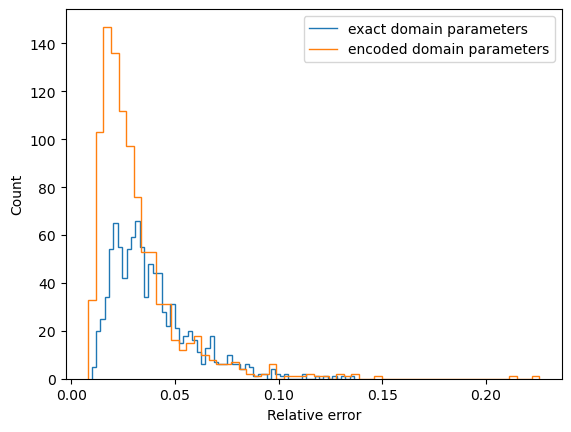}
     \caption{Histograms of relative errors on test set for model with ROM encoding dimension 30, trained on the exact ellipse parameters (blue) and model with ROM encoding dimension 30, trained on 20-dimensional domain encoding (orange). When using exact parametrization, $15$ samples have the relative error larger than $0.1$ and $213$ larger than $0.05$. In case of using calculated domain parametrization, $15$ samples have the relative error larger than $0.1$ and $121$ larger than $0.05$. There are more error outliers when using exact ellipse parametrization. }
     \label{elipseDist}
 \end{figure}

 In order to   compare the exact and autoencoder domain parametrization, we investigate error distributions for models with the best performance according to Table \ref{table1}. Histograms of relative errors are in Figure \ref{elipseDist}. The results indicate that domain encoding improves model robustness, as reflected by a reduced presence of outliers.

\subsection{Domains with varying number of parameters}
We experiment again with the advection-diffusion equation on parametric domains. In this setting, the domains vary in complexity by depending on a different number of parameters, which makes them challenging to parametrize using classical techniques. Specifically, we consider square domains that contain between one and four circular holes. More precisely, we consider the following boundary value problem:
\begin{equation}
    \begin{cases}
      -\Delta u + 0.1 \mathbf{c} \cdot \nabla u = 1 & \text{, on $\langle0, 1\rangle^2$ without circles}\\
      u = 2 &\text{, on circles' boundaries} \\
      u = 0 &\text{, for $y=0$ and $y=1$}\\
      \frac{\partial u}{\partial \nu} = \beta & \text{, for $x=0$} \\
      \frac{\partial u}{\partial \nu} = -\beta & \text{, for $x=1$}      
    \end{cases}       
\end{equation}
Parameters are:

\begin{itemize}
    \item  transport angle $\phi \in   [0, 2 \pi]$ which determines vector $\mathbf{c}=(cos(\phi), sin(\phi))$ 
    \item $m \in [1,4]$ number of circular holes
    \item  centers of holes $(x_i, y_i)$, $i = 1, \ldots m$
    \item circles' radii $r_i$, $i=1,\ldots m$
    \item boundary flux $\beta \in [1, 10]$
\end{itemize}
In total, we again have two PDE parameters: transport angle $\phi$ and boundary flux $\beta$, but there are between three and twelve geometric parameters, namely, the radius and center coordinates for each circular hole.

\begin{figure}[H]
     \centering
     \includegraphics[scale=0.33]{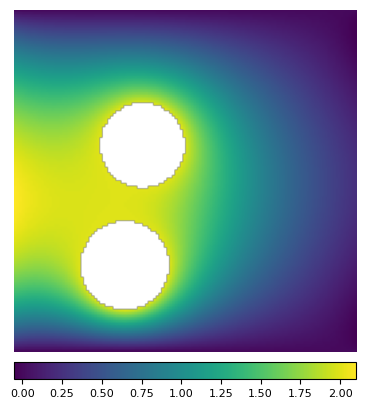}\hspace{3ex}
        \includegraphics[scale=0.33]{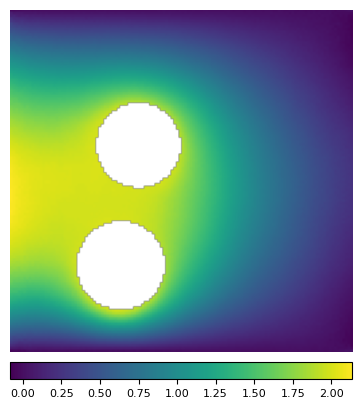}
     \includegraphics[scale=0.33]{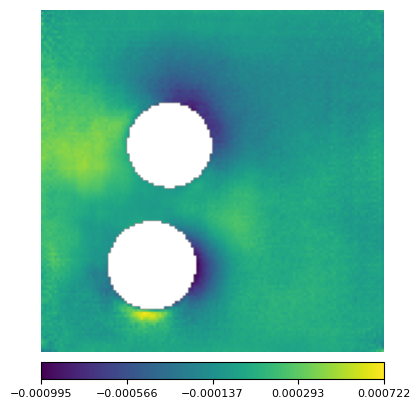}
     \caption{Exact solutions (left), model prediction (middle), relative error (right) for model trained on $1000$ random samples for parameters $\beta = 1.671$, $\phi=0.182$,  with centers $(0.322, 0.379)$ and $(0.374, 0.736)$, and radii $0.130$ and $0.124$, respectively. }
    \label{circles}
\end{figure}

\begin{figure}[H]
     \centering
     \includegraphics[scale=0.33]{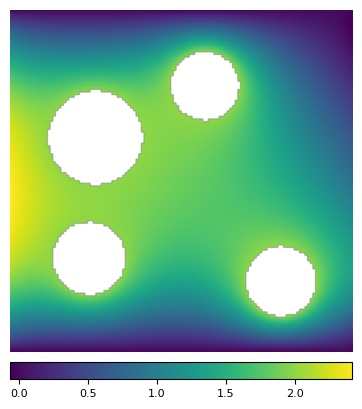}\hspace{3ex}
        \includegraphics[scale=0.33]{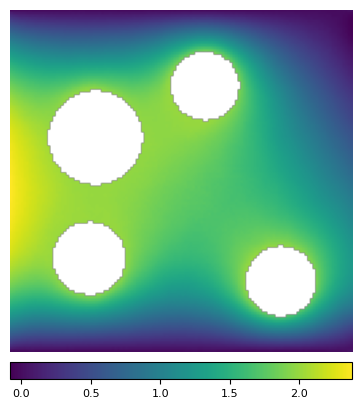}
     \includegraphics[scale=0.33]{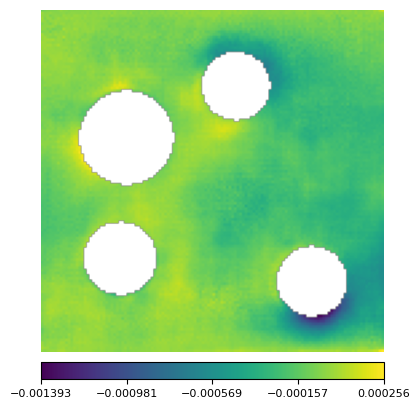}
     \caption{Exact solutions (left), model prediction (middle), relative error (right) for model trained on $1000$ random samples for parameters $\beta = 3.204$, $\phi=2.893$,  with centers $(0.229, 0.729)$, $(0.569, 0.222)$, $(0.248, 0.374)$ and $(0.792, 0.796)$, and radii $0.107$, $0.101$,$0.139$,$0.103$
, respectively. }
    \label{circles2}
\end{figure}

To generate domains, first, we randomly choose the number of circles $m$. Then we randomly choose the centers of the circles so that the distance between the centers and the outer boundary is at least $0.1$. Finally, radii are chosen randomly so that the distance between all circles is at least $0.1$ and $r \in [0.05, 0.2]$. This procedure ensures that circles do not overlap and that all circles are of reasonable size. We have chosen the same encoder-decoder architectures as in the previous example, and after training solution autoencoders, a grid search was performed to obtain $\Phi_S$ hyperparameters.
Table \ref{table3} shows the average relative error on the test set for different combinations of domain and solution latent dimensions, while Figures \ref{circles} \ref{circles2} show examples of solutions corresponding to a geometric configurations from the test set.
In this case, constructing classical ROM is unfeasible because of different topologies and placements of the circles.

\begin{center}
\begin{table}[H]
\begin{adjustbox}{ center}
\begin{tabular}{ |c |c |c |c| c| }
 \hline
  &  \thead{ test error, \\ ROM encoding \\ dimension 15, \\ domain encoding \\ dimension 20}
& \thead{test error, \\ ROM encoding \\ dimension 15, \\ domain encoding \\ dimension 30} & 
\thead{test error, \\ ROM encoding \\ dimension 30, \\ domain encoding \\ dimension 20}  
&\thead{test error, \\ ROM encoding \\ dimension 30, \\ domain encoding \\ dimension 30} \\
 \hline 
    \thead{calculated \\
    domain parameters}                   &  
    0.0252          &  0.0242        & 0.0244    & 0.0234   \\
    \hline
    
\end{tabular}
\end{adjustbox}
    \caption{Test relative errors obtain using different ROM and the domain encoding dimensions. In this case, domains have varying numbers of parameters. We use only geometric parameters obtained by the autoencoder. We report on the best model for each ROM encoding dimension and domain encoding dimension. A smaller domain encoding dimension is beneficial in this case. }
    \label{table3}
\end{table}
\end{center}

\textbf{Robustness to Domain Deformations.} 
To test the robustness of the model, we evaluated DL-ROMs trained on domains with circular holes using test domains in which the circles were slightly deformed into ellipses. Specifically, each circle was transformed by reducing its minor semiaxis to $(1 - p)$ times the major semiaxis, where $p$ is a deformation factor. To introduce variability, even-numbered ellipses were aligned along the $x$-axis, and odd-numbered ones were aligned along the $y$-axis. This setup simulates realistic domain perturbations, such as geometric noise or measurement uncertainty.

The trained models were tested on these deformed domains, and their predictions were compared with finite element solutions. Table~\ref{table4} reports the average test errors for increasing levels of deformation. As expected, the error increases as the deformation grows. However, the model does not break down and continues to produce reasonable approximations even when tested on domains with types of geometries not seen during training. This suggests that the method generalizes well to moderate domain perturbations.

\begin{figure}[h]
     \centering
     \includegraphics[scale=0.33]{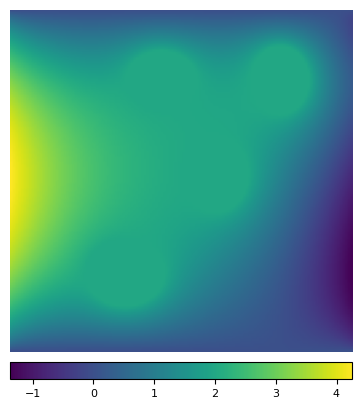} 
        \includegraphics[scale=0.33]{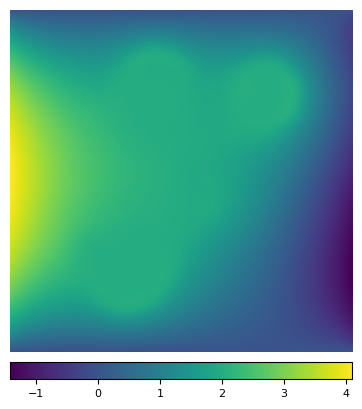}
     \includegraphics[scale=0.33]{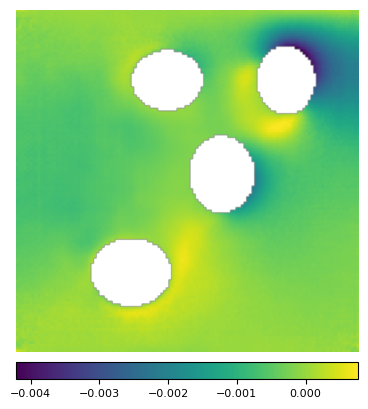}
  \caption{The exact solution (left) not multiplied by characteristic function, the model prediction (middle) not multiplied by a characteristic function, relative error (right) multiplied by characteristic function for slightly deformed circles with parameters $\beta = 9.643$, $\phi = 4.129$. Ellipses have centers, $(0.601, 0.0481)$, $(0.441, 0.0.206)$, $(0.0.789,0.0.206)$, and $(0.334, 0.770)$, and radii of non-deformed circles $0.113$, $0.105$, $0.101$, and $0.118$, respectively. Small half-axes are $15\%$ smaller than large half-axes for each ellipse. While relative error is highest near domain boundaries due to unseen shapes during training, the model maintains low error elsewhere, showing robustness to domain perturbation.
  }
    \label{circlesSensitivity}
\end{figure}

\begin{figure}[H]
    \centering
    \includegraphics[width=0.4\linewidth]{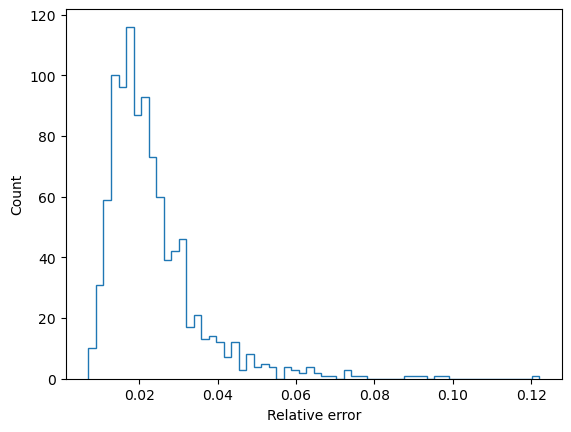}    
    \caption{Histogram of relative errors on test set for model with ROM encoding dimension 30 and domain encoding dimension 30. Only $38$ out of $1000$ samples have relative error larger than $0.05$.}
    \label{circlesDist}
\end{figure}

Figure~\ref{circlesSensitivity} shows one such out-of-distribution example. The relative error is highest near the boundaries of the ellipses, where the geometry most significantly deviates from the training distribution. Still, the model captures the global solution behavior accurately. The histogram of test errors in Figure \ref{circlesDist} confirms that even under moderate distribution shifts, no significant outliers appear, indicating that the predictions remain stable and that the method is robust.

\begin{center}
\begin{table}[H]
\begin{adjustbox}{ center}
\begin{tabular}{ |c |c |c |c| c| }
 \hline
  &  \thead{ test error, \\ ROM encoding \\ dimension 15, \\ domain encoding \\ dimension 20}
& \thead{test error, \\ ROM encoding \\ dimension 15, \\ domain encoding \\ dimension 30} & 
\thead{test error, \\ ROM encoding \\ dimension 30, \\ domain encoding \\ dimension 20}  
&\thead{test error, \\ ROM encoding \\ dimension 30, \\ domain encoding \\ dimension 30} \\
 \hline 
    \thead{$p = 0$}                   &   0.0252
              & 0.0242         & 0.0244       &   0.0234\\
    \hline
        \thead{$p = 0.05$}                   & 0.0260
              &  0.0252    &  0.0252    & 0.0254    \\
    \hline
        \thead{$p = 0.1$}                   &  0.0280
               &   0.0275       &  0.0278      & 0.0294    \\
    \hline
        \thead{$p = 0.15$}                   &  0.0309
               &    0.0310     &   0.0319     &   0.0356\\
    \hline
    
\end{tabular}
\end{adjustbox}
    \caption{Relative test errors for DL-ROMs trained on domains with circular holes ($p = 0$), evaluated on domains with increasingly deformed elliptical holes. Here, $1 - p$ denotes the ratio between the smaller and larger semi-axis of the ellipse. We report results for the best model for each ROM and domain encoding dimension. Performance degrades gradually with increasing deformation, but remains acceptable, demonstrating robustness to moderate boundary changes.}
    \label{table4}
\end{table}
\end{center}

\subsection{Navier-Stokes equations}

In this section, we evaluate the DL-DA-ROM framework on the Navier–Stokes equations, a significantly more complex and challenging problem. The Navier–Stokes equations govern fluid flow and include significant nonlinearities, providing a challenging and realistic scenario to test the effectiveness and robustness of our methodology in practical applications.

The considered scenario involves fluid flow through a cavity with a geometrically complex boundary, motivated by the simulations of blood flow through vessels with aneurysms. Specifically, the domain is characterized as a channel with a variable width, representing an aneurysm-like enlargement, and further complicated by the presence of a circular obstacle to introduce additional geometric complexity. We vary the Reynolds number between 1 and 100.

\begin{figure}[H]
     \centering
     \includegraphics[scale=0.5]{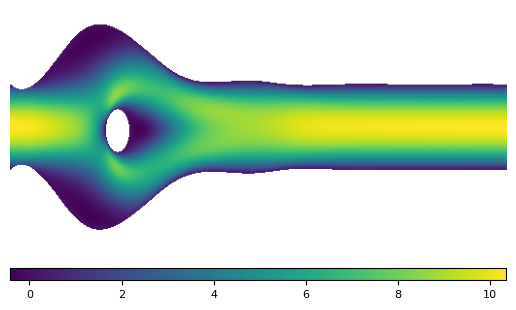}\hspace{3ex}
        \includegraphics[scale=0.5]{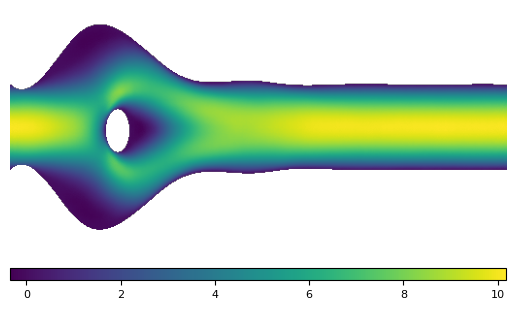}
     \includegraphics[scale=0.5]{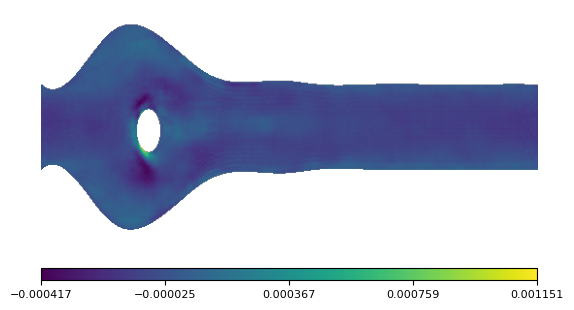} 
     \caption{The exact solutions (top left), model prediction (top right) and relative error (bottom) for model trained on $1000$ random parameter samples, with $20$-dimensional domain parametrization obtain using $a = 10.34$ and $Re = 30.1$. }
    \label{NS}
\end{figure}

\begin{figure}
    \centering
    {\includegraphics[scale = 0.5]{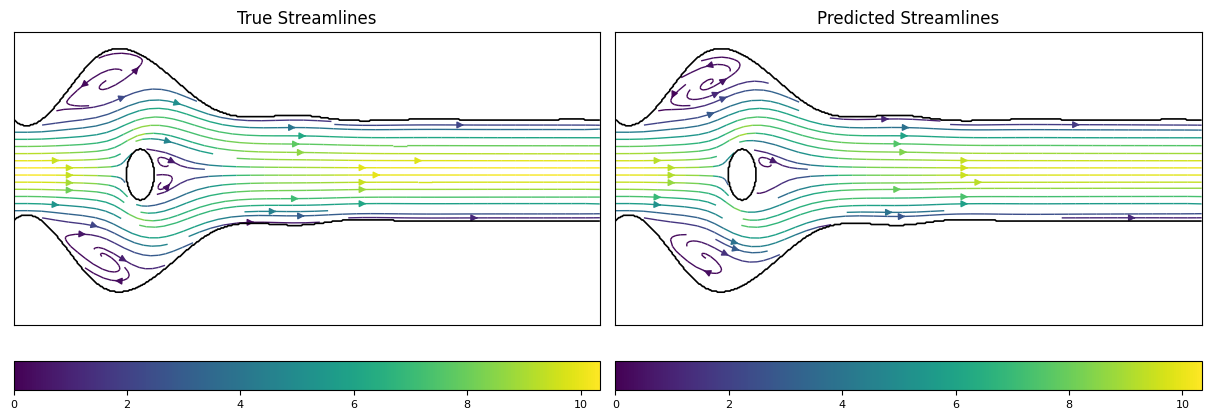}}
    \caption{Correct streamlines (left) and predicted (right) for $a = 10.34$ and $Re = 30.1$. }
    \label{streamlinePlot}
\end{figure}

Furthermore, to assess the robustness of our model, we deliberately test parameter configurations that lie outside the original training distribution. This allows us to verify that our model can handle moderate extrapolation scenario, ensuring stability and reasonable accuracy for previously unseen, but not drastically different, geometries and flow conditions. Precisely, we consider the following boundary value problem:
    \begin{equation*}
    \begin{cases}
        (\boldsymbol{u} \cdot \nabla)\boldsymbol{u} - \frac{1}{Re}\Delta \boldsymbol{u} = \nabla p &\text{in $\Omega(\lambda)$,} \\
      \boldsymbol{u} = 0 &\text{on the obstacle and top and bottom walls,} \\
      \boldsymbol{u}= (a p(y), 0) &\text{on the inflow, where p is a parabolic function,}\\
      \boldsymbol{\sigma n} = 0
      &\text{on the outflow},  
    \end{cases}       
\end{equation*}
where the domain $\Omega(\lambda))$ is a cavity without a circle. The length of the domain is $L=4.5$ and $\sigma$ is the fluid stress tensor. Top side is deformed as $R + h e^{- w(x-1.25)^2}+ s(x)$, where the spline $s$ is given by $9$ points $(0,0), (0.25, 0), (x_{s1}, y_{s1}), (x_{s2}, 0), (2.2,0), (x, h e^{(-w( -L/2)^2)} - he^{-w(x-L/2)^2})$ for $x = 3.0, 3.5, 4.0, 4.5$ and the bottom is deformed as $-he^{ w(x-1.25)^2} - s(x)$.
The parameters are: 
\begin{itemize}
    \item undeformed pipe width $R \in [0.1, 0.2]$ 
    \item obstacle radius $rCirc \in [0.1,0.15]$ 
    \item $x_{s1} \in [0.5, 1 ]$ 
    \item $y_{s1} \in [0.15, 0.25]$
    \item $x_{s2} \in [1.5, 2]$ 
    \item obstacle center $(x_0, y_0) \in [0.5,1.5] \times [0.05, 0.1]$ 
    \item inflow amplitude $a \in [5,15]$
    \item $w \in [0.05, 0.2]$
    \item $h \in [0.15, 0.2]$ 
    \item $Re \in [1, 100]$
\end{itemize}
In total, we have two PDE parameters: inflow amplitude $a$ and Reynolds number $Re$ and $9$ geometric parameters: pipe width $R$, obstacle radius $rCirc$, obstacle center $(x_0, y_0)$, spline points $x_{s1}, y_{s1},x_{s2}$ and $w$ and $h$ determining the shape of the cavity.
 The solutions are interpolated on the $256 \times 512$ grid on $[-0.5, 0.725] \times [0, 4.5]$. We train the model to predict the velocity in $x$ and $y$-direction. The autoencoder is trained on the velocity components in both the $x$ and $y$-directions, while all reported relative errors correspond to the velocity magnitude. A comparison of model performance using calculated versus exact domain parameters is presented in Table~\ref{table_ns}, with the representative results shown in Figure~\ref{NS} and streamlines are shown in Figure \ref{streamlinePlot}. We emphasize that despite the increased complexity of both the governing equations and the domain geometry compared to the advection-diffusion case, the average error only slightly increases and remains well within acceptable limits for reduced-order modeling.
 
 The geometry models a simplified fusiform aneurysm with a localized symmetric enlargement of the vessel and an idealized intraluminal inclusion, represented as a circular obstacle, intended to mimic features such as thrombus or calcified deposits. The parabolic inflow drives flow through the channel, and across the simulated parameter range the flow remains steady and attached to the gradually curved walls. The dominant visible features are a jet that widens as it enters the enlargement and contracts again downstream, and a low-velocity wake behind the obstacle. This class of geometries is relevant to hemodynamic simulations of fusiform aneurysms, where the interaction between the enlarged vessel section and intraluminal features influences the local flow field.

 As in the previous examples, a standard ROM based on mapping all domains to a single reference configuration is not directly applicable here. Although the wall deformations in this test are moderate, the simultaneous variation of the channel shape and the circular obstacle would still require a problem-specific reference mapping and associated mesh-handling procedure in order to construct a classical POD-Galerkin ROM. For this reason, in this example we do not include a comparison with such methods and instead focus on data-driven baselines that can be applied directly to the interpolated fields.

\begin{center}
\begin{table}[ht] 
\begin{adjustbox}{   center}
\begin{tabular}{ |c |c |c |c|c| }
 \hline
  
 & \thead{test error, \\ ROM encoding \\dimension 15, \\ domain encoding \\ dimension 20} 
 & \thead{test error, \\ ROM encoding \\ dimension 15, \\ domain  encoding \\ dimension 30} 
  &  \thead{test error, \\ ROM encoding \\ dimension 30, \\ domain encoding \\ dimension 20}
& \thead{test error, \\ ROM encoding \\ dimension 30, \\  domain encoding \\ dimension 30}   
\\
\hline
  \thead{exact domain \\ parameters only} & 0.0412 
            &   0.0412   & 0.0335   & 0.0335            
                             \\
   \hline
    \thead{calculated \\
    domain parameters} &  0.0361              
                      & 0.0348       &  0.0379       & 0.0333    \\
    \hline
\end{tabular}
\end{adjustbox}
    \caption{Test relative errors obtained using different ROM and the domain encoding dimensions. We differentiate between cases where geometric parameters are exactly known (1st row) and when we use only geometric parameters obtained by the autoencoder (2nd row). We can see that the use of the autoencoder domain parameterization does not negatively affect DL-ROM performance. We report on the best model for each ROM encoding dimension and domain encoding dimension. Most Models trained on calculated domain parameters perform better than the models trained using only exact domain parameters. Even though results are similar, larger fully-connected neural network is needed in case of calculated domain parameters (4 layers with 2048 neurons vs 4 layers with 512 neurons).
  }
    \label{table_ns}
\end{table}
\end{center}

We further evaluate the sensitivity of the model to domain perturbations by testing parameter configurations that fall outside the range of the original training set. Specifically, we consider domains with wider cavities ($w \in [0.2, 0.25]$), using the model with ROM and domain encoding dimensions both set to $30$. Additionally, we increase the Reynolds number slightly beyond the training range, with $Re \in [100, 110]$, while other parameters remain within the original bounds. A total of $100$ such samples were generated, resulting in an average relative error of approximately $0.065$.

\begin{figure}[H]
    \centering
    \includegraphics[width=0.4\linewidth]{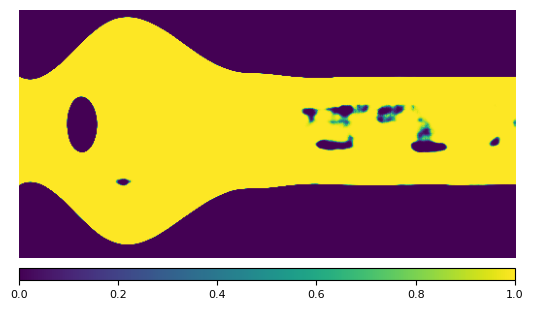}
        \includegraphics[width=0.4\linewidth]{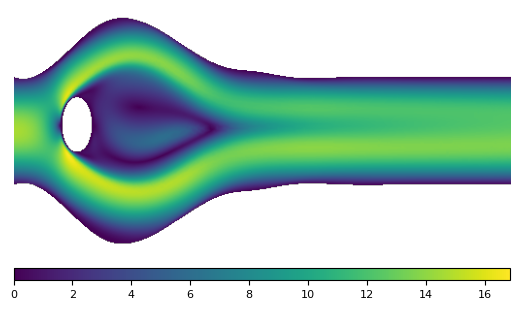}
    \includegraphics[width=0.4\linewidth]{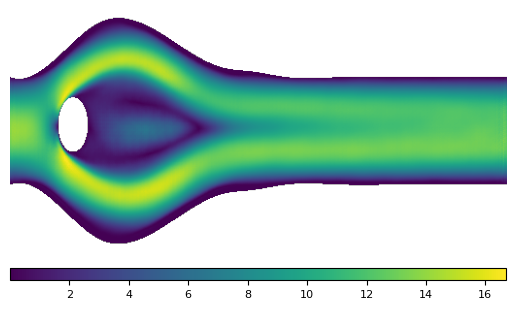}
    \includegraphics[width=0.4\linewidth]{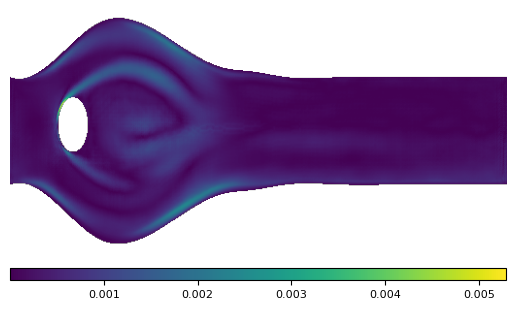}
    \caption{Out of distribution sample with the highest relative error of $0.138$. Domain reconstruction (top left), correct solution (top right), prediction (bottom left), relative error (bottom right).  Domain autoencoder is not able to properly reconstruct the domain, producing artificial holes in the domain. This implies that calculated domain parametrization is not appropriate for this domain. However, the model still finds a reasonable solution. Even though the domain reconstruction is incorrect, the solution autoencoder decoder is constrained to the manifold of solutions seen during training. It therefore produces the closest plausible solution within that learned distribution, which happens to still be physically reasonable, thus explaining the moderate rather than catastrophic error.}
    \label{wrongReconstruction}
\end{figure}

\begin{figure}[H]
    \centering
    \includegraphics[width=0.4\linewidth]{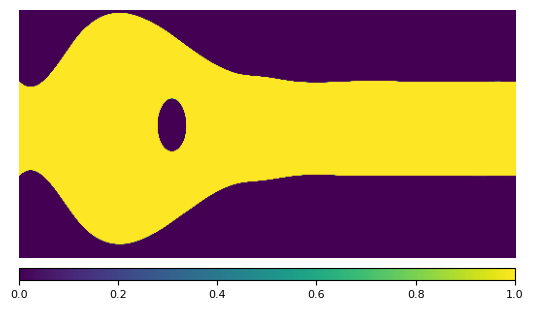}
        \includegraphics[width=0.4\linewidth]{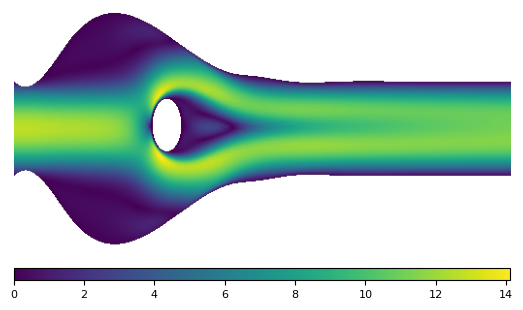}
    \includegraphics[width=0.4\linewidth]{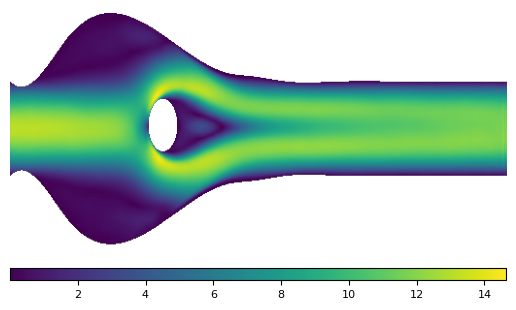}
    \includegraphics[width=0.4\linewidth]{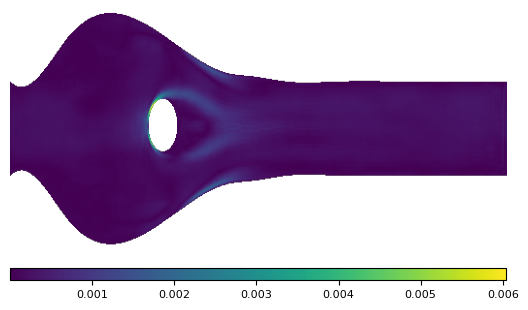}
    \caption{Typical out of distribution sample. Domain reconstruction (top left), correct solution (top right), prediction (bottom left), relative error (bottom right). In this case, the domain autoencoder is able to reconstruct the domain, thus able to extract more meaningful representation.}
    \label{okReconstruction}
\end{figure}

We closely examine outlier cases and find that one source of increased error occurs when the domain autoencoder fails to accurately reconstruct the domain geometry (see Figure~\ref{wrongReconstruction}). Nevertheless, even in such cases, the method provides a reasonable approximation with a relative error of $0.138$. Conversely, Figure~\ref{okReconstruction} illustrates a case where the domain encoding performed well, resulting in a relative error of $0.087$. These examples indicate that, despite some reconstruction issues, the model remains robust and produces acceptable solutions even under moderate parameter and geometry shifts.

\begin{figure}[H]
    \centering
    \includegraphics[width=0.45\linewidth]{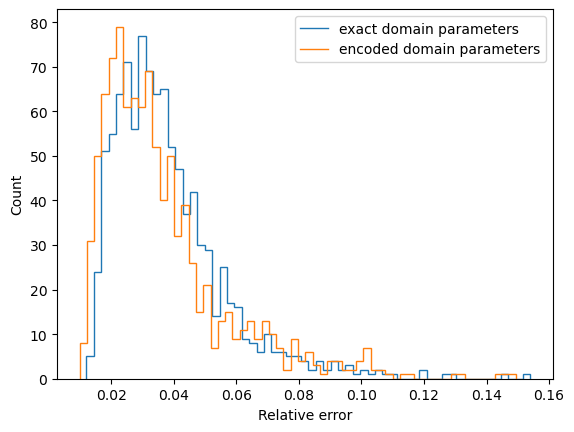}     
    \caption{Histogram of relative errors on test set for the model with ROM encoding dimension 30 trained with exact domain parameters (blue), and trained with domain encoding dimension 30 (orange).  }
    \label{circlesDistNS}
\end{figure}

The distribution of the relative error on the test set for the best model reported in Table \ref{table_ns} is shown in Figure \ref{circlesDistNS}. We can see that the model predictions are stable across the test dataset. 

\subsection{Comparison with POD-based methods}
In this section, we compare our model against POD-based methods. In all experiments, classical ROM with POD and mesh mapping is infeasible, as mapping solutions to a reference domain is numerically intractable, either due to the high variability of the domains involved or differences in topology, as discussed for each PDE separately. POD projection and POD-DL-ROM are instead implemented using POD representations of interpolated solutions. The POD decay is moderate across all examples: between 30 and 50 modes are required to achieve the same reconstruction error as DL-DA-ROM. 

\begin{figure}[!h]
    \centering
    \includegraphics[width=0.33\linewidth]{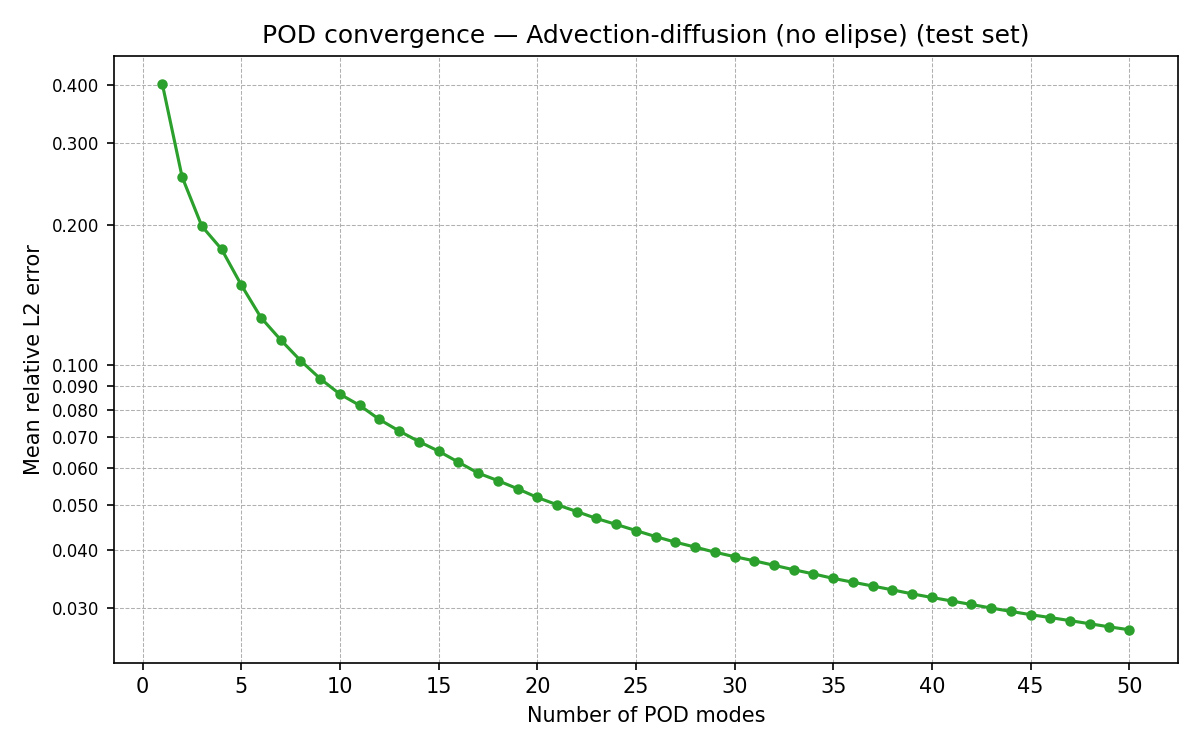}
    \includegraphics[width=0.33\linewidth]{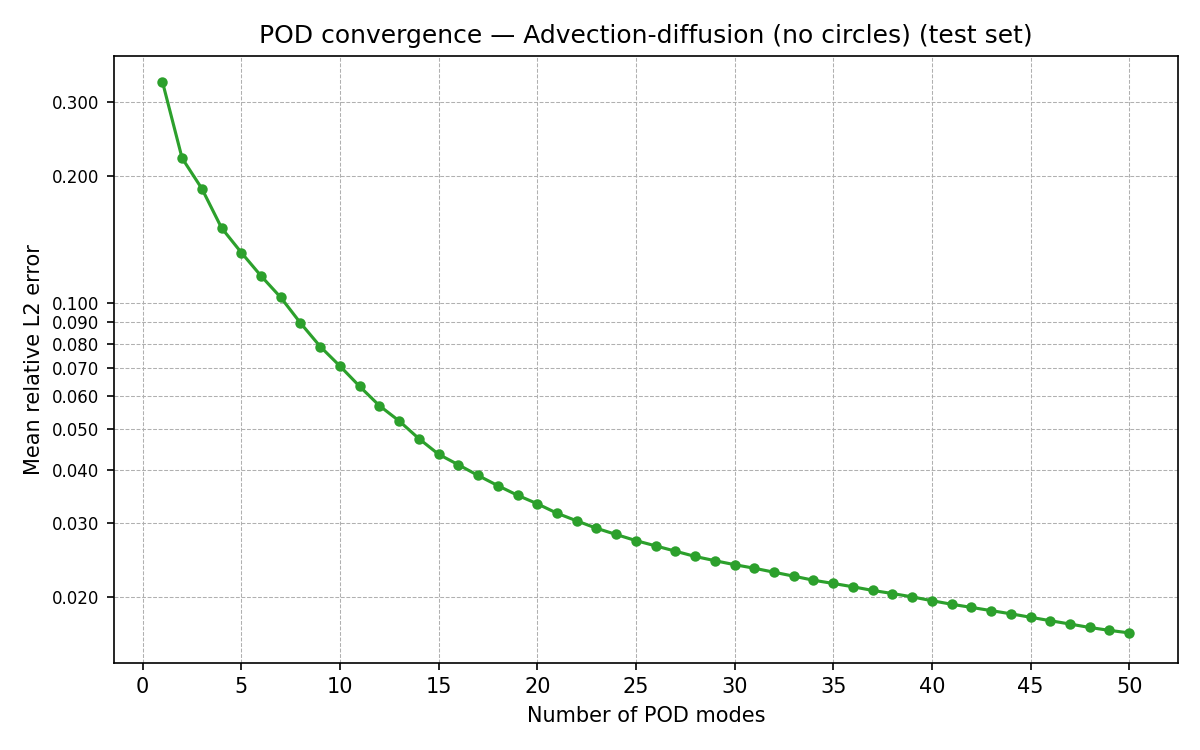}   
    \includegraphics[width=0.33\linewidth]{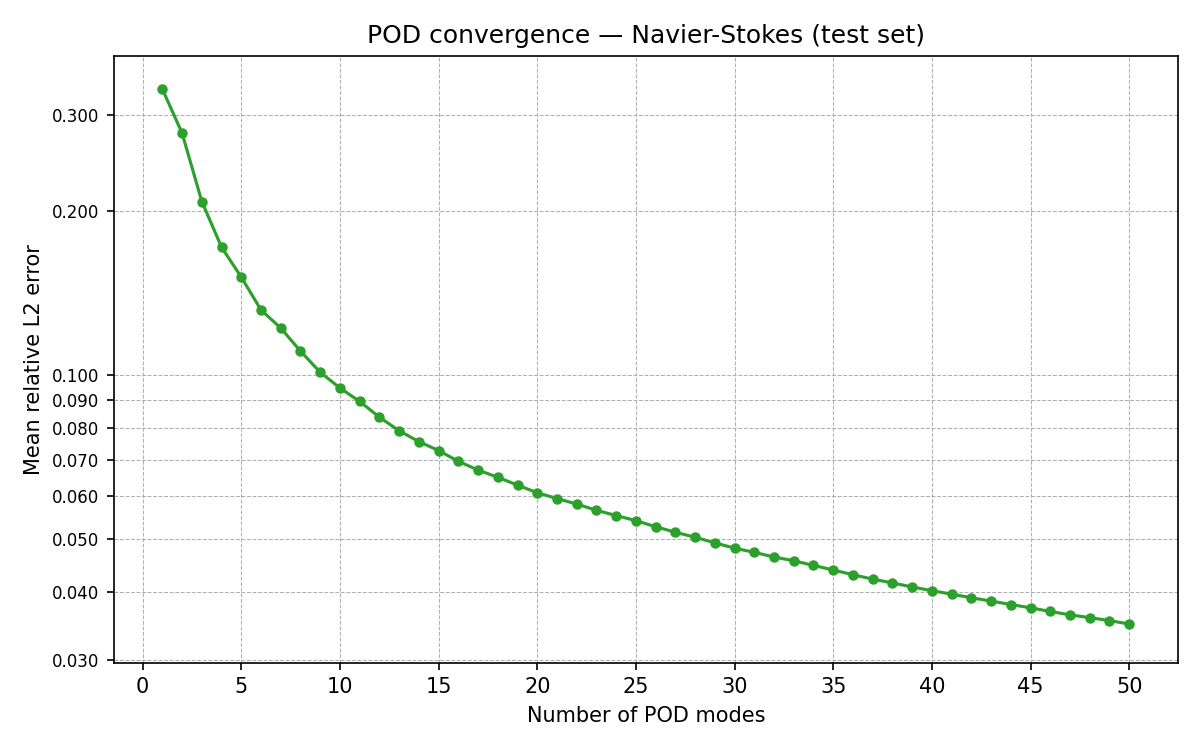}
    \caption{POD decay of interpolated solutions on the test set. The number of modes required to match DL-DA-ROM reconstruction error is 41 for advection-diffusion on a square domain without an ellipse, 31 for advection-diffusion on domains without one to four circles, and 50 for the Navier-Stokes case.}
    \label{podDecay}
\end{figure}

While such a POD basis could in principle support data-driven methods, it cannot sustain a classical POD-Galerkin projection in our setting. The basis is constructed on a fixed interpolation grid that is disconnected from the original FEM spaces on varying domains, and no consistent reference configuration exists across samples. Nevertheless, these POD modes can be leveraged by POD-based data-driven methods such as POD-DL-ROM, as illustrated in Figure \ref{podDecay}.

To ensure a fair comparison with POD-DL-ROM, we performed a hyperparameter grid search to identify the best configuration. As with DL-DA-ROM, the latent representation of the characteristic function, obtained from the same domain autoencoder, is used as the domain parametrization. The selected hyperparameters are reported in Table \ref{poddlromgrid}.

 \begin{center} 
\begin{table}[H] \centering
\begin{tabular}{ |c |c |}
 \hline
  
    Hyperparameter  & Values \\ \hline
  Training set size & 1000 \\ \hline 
  Number of epochs & 1000 \\ \hline 
  Number of POD modes & 64, 256, 400 \\ \hline
  Number of layers & 4 \\ \hline 
  Number of channels per CNN layer & 16, 32 , 64 \\ \hline
  Neurons per layer for $\Phi_S$ & 512, 1024, 2048 \\ \hline
  Maximum learning rate & 1e-4 \\ \hline
  Batch size & $50$ \\ 
 \hline 
\end{tabular}
\caption{Hyperparameters  grid search for POD-DL-ROM. Encoder, decoder and $\Phi_S$ have $4$ layers, where $\Phi_S$ is fully-connected neural network that maps parameters to solution latent representation, as in DL-DA-ROM. Encoder and decoder are convolutional. POD coefficients are reshaped into square grid prior to input to the encoder.   }
    \label{poddlromgrid}
\end{table}
\end{center}

We report the relative test error for POD-DL-ROM in Table \ref{poddlcomparison} for the model with the lowest relative validation error.
\begin{center}
    \begin{table}[H] \centering 
        \centering
        \begin{tabular}{|c|c|c|c|} \hline 
             & Advection-diffusion (ellipse) & Advection-Diffusion (circles) & Navier-Stokes \\ \hline 
             POD-DL-ROM& 0.2297 &0.0932 & 0.1245 \\ \hline 
            DL-DA-ROM & \textbf{0.0310} & \textbf{0.0234} & \textbf{0.0333}  \\  \hline 
        \end{tabular}
        \caption{Comparison of POD-DL-ROM and DL-DA-ROM. DL-DA-ROM outperforms POD-DL-ROM across all three test cases. We attribute this to the fact that DL-DA-ROM operates directly on the interpolated field, where neighbouring pixels correspond to neighbouring points in physical space, which is a structure that matches what convolutional networks are designed to exploit. POD coefficients, in contrast, represent the solution as projections onto global basis functions, and this representation appears harder to learn from when domain geometry varies significantly across samples. }
        \label{poddlcomparison}
    \end{table}
\end{center}
We conclude that, for the class of parameter-dependent domain problems considered here, DL-DA-ROM is not only applicable in settings where classical POD-Galerkin ROM is impractical, but also consistently outperforms a carefully tuned POD-DL-ROM baseline in terms of accuracy.

\section{Conclusion}
We presented a method for building reduced-order models for PDEs on a parameter-dependent domains, including scenarios where the domain geometry is not explicitly given but instead defined through images or measurements. Unlike many existing ROM techniques that require consistent meshing or graph-based representations, our approach directly utilizes bitmap images of domain characteristic functions. This makes it particularly suitable for practical applications where mesh generation is costly, unreliable, or unavailable—such as simulations derived from imaging data or experimental measurements.

Our results demonstrate that the learned domain encoding is not only compact but also meaningful; in cases with known geometric parameters, our encoding achieved similar accuracy to models using exact parameterizations. Furthermore, our method effectively handled significant variations in domain shape, changes in the number of features (such as holes), and moderate boundary perturbations without losing stability.
Although the use of domain parametrization obtained via the CNN autoencoder successfully enables mapping from physical parameters to the latent solution space, we observed that this configuration generally requires a larger network capacity to achieve comparable accuracy. This behavior likely stems from a partial misalignment between the domain parametrization and the reduced-order latent space, which increases the effective nonlinearity of the mapping.
Future work will focus on improving the consistency between the autoencoder latent representation and the parameter space, for example through joint or coupled training strategies or latent regularization terms. Such approaches could reduce the required network size and improve generalization while preserving the interpretability and robustness of the proposed DL-DA-ROM framework.
Conceptually, our approach extends naturally to three-dimensional problems. However, from a practical perspective, computational efficiency remains an important challenge due to the increased size of the input grids. Improving computational efficiency for three-dimensional implementations is currently a focus of our ongoing work.

Additionally, the learned domain representations may offer utility beyond DL-ROM frameworks, potentially benefiting other situations requiring domain parameterization where explicit parameters are unavailable or in different machine learning-based model reduction contexts.

In short, our contribution can be summarized as follows:
\begin{itemize}
    \item Introduced a CNN autoencoder-based domain parameterization that handles challenging scenarios including topological changes (varying numbers of holes/features) and domains defined solely through imaging or measurement data.
    \item Proposed a DL-DA-ROM framework that directly processes domains given on discrete rectangular grid, eliminating the need for mesh generation or explicit geometric parameterization.
    \item Validated the framework on multiple PDEs across domains with varying geometry and topology, showing accuracy comparable to models using exact geometric parameterizations.
\end{itemize}

\section{Acknowledgments}
This research was carried out using the advanced computing service provided by the University of Zagreb University Computing Centre - SRCE. This research was supported by the Croatian Science Foundation under the project number IP-2022-10-2962 (I.M., B.M. and D.V.). M.B. has been partly supported by the National Science Foundation via grants NSF DMS-2208219 and NSF DMS-2205695. BM was supported by the European Union – NextGenerationEU through the National Recovery and
Resilience Plan 2021-2026. Institutional grant of University of Zagreb Faculty of Science IK IA 1.1.3. Impact4Math.

We thank the anonymous referees for their careful reading and valuable comments, which helped improve the manuscript.

\bibliographystyle{unsrt}  
\bibliography{references}

\end{document}